\documentclass[reqno,twoside,12pt]{amsart}

\hoffset - 1.5cm

\usepackage{graphicx}
\usepackage{amstext}
\usepackage[T1]{fontenc}
\usepackage[cp1250]{inputenc}
\usepackage{pstricks,pst-node,pst-text,pst-3d}
\usepackage{color}
\usepackage{amsmath}
\usepackage{amsthm}
\usepackage{amssymb}
\usepackage[mathscr]{euscript}
\usepackage{bbm}
\usepackage{enumerate}
\usepackage{arydshln}
\usepackage{MnSymbol}
\newtheorem{Remark}{Remark}[section]
\newtheorem{Corollary}[Remark]{Corollary}
\newtheorem{Definition}[Remark]{Definition}
\newtheorem{Example}[Remark]{Example}
\newtheorem{Fact}[Remark]{Fact}
\newtheorem{Lemma}[Remark]{Lemma}
\newtheorem{Proposition}[Remark]{Proposition}
\newtheorem{Theorem}[Remark]{Theorem}

\newcommand{\ba}{\begin{array}}
\newcommand{\bc}{\begin{center}}
\newcommand{\bd}{\begin{description}}
\newcommand{\bdm}{\begin{displaymath}}
\newcommand{\be}{\begin{enumerate}}
\newcommand{\beq}{\begin{equation}}
\newcommand{\bdf}{\begin{Definition}}
\newcommand{\bex}{\begin{Example}}
\newcommand{\bft}{\begin{Fact}}
\newcommand{\bl}{\begin{Lemma}}
\newcommand{\bp}{\begin{Proposition}}
\newcommand{\br}{\begin{Remark}}
\newcommand{\bt}{\begin{Theorem}}
\newcommand{\bco}{\begin{Corollary}}
\newcommand{\bhy}{\begin{Hypothesis}}
\newcommand{\ea}{\end{array}}
\newcommand{\ec}{\end{center}}
\newcommand{\ed}{\end{description}}
\newcommand{\edm}{\end{displaymath}}
\newcommand{\ee}{\end{enumerate}}
\newcommand{\eeq}{\end{equation}}
\newcommand{\edf}{\end{Definition}}
\newcommand{\eex}{\end{Example}}
\newcommand{\eft}{\end{Fact}}
\newcommand{\el}{\end{Lemma}}
\newcommand{\ep}{\end{Proposition}}
\newcommand{\er}{\end{Remark}}
\newcommand{\et}{\end{Theorem}}
\newcommand{\eco}{\end{Corollary}}
\newcommand{\ehy}{\end{Hypothesis}}

\newcommand{\bN}{\mathbb{N}}

\newcommand{\bR}{\mathbb{R}}
\newcommand{\bS}{\mathbb{S}}

\newcommand{\bV}{\mathbb{V}}

\newcommand{\bZ}{\mathbb{Z}}

\newcommand{\cB}{\mathcal{B}}
\newcommand{\cC}{\mathcal{C}}

\newcommand{\cH}{\mathcal{H}}

\newcommand{\cN}{\mathcal{N}}

\newcommand{\cT}{\mathcal{T}}
\newcommand{\cU}{\mathcal{U}}

\newcommand{\Comega}{C_{2\pi}(\Omega)}

\newcommand{\diag}{\mathrm{ diag \;}}
\numberwithin{equation}{section} \errorcontextlines=0

\newcommand{\sign}{\mathrm{ sign \;}}

\newcommand{\card}{\mathrm{card}\:}
\newcommand{\morse}{\mathrm{m^-}}

\newcommand{\cl}{\mathrm{cl}}

\newcommand{\sone}{S^1}

\newcommand{\ds}{\displaystyle}
\newcommand{\nt}{\noindent}

\newcommand{\ib}{\mathrm{i}_{\mathrm{B}}}

\newcommand{\bif}{\mathcal{BIF}}

\usepackage{todonotes}
\usepackage{setspace}

\makeatletter
\@namedef{subjclassname@2020}{%
  \textup{2020} Mathematics Subject Classification}
\makeatother

\begin{document}

\title[Periodic solutions of Hamiltonian systems]{Generalization of Lyapunov Center Theorem for Hamiltonian systems via normal forms theory}

\author{Anna Go{\l}\c{e}biewska$^{1)}$}
\address{$^{1), 2)}$Faculty of Mathematics and Computer Science \\ Nicolaus Copernicus University in Toru\'n\\
PL-87-100 Toru\'{n} \\ ul. Chopina $12 \slash 18$ \\
Poland}

\author{S{\l}awomir Rybicki$^{2)}$}

\email{aniar@mat.umk.pl (A. Go{\l}\c{e}biewska)}
\email{rybicki@mat.umk.pl (S. Rybicki)}

\date{\today}


\begin{abstract}
In this article we formulate and prove  sufficient   conditions for the existence of trajectories of nonstationary periodic solutions of autonomous Hamiltonian systems in a neighbourhood of equilibria. It is worth pointing out that assumptions of some well-known theorems imply that of our main results. We obtain our results with the use of the theory of normal forms for Hamiltonian matrices and global bifurcation theory for autonomous Hamiltonian systems.
\end{abstract}
%%%%%%%%%%%%%%%%%%%%%%%%%%%%%%%%%%%%%%%%%%%%%%%%%%%%%%%%%%%%%%%%%%%%%%%%%%%%%%%%%%%%%%%%%%%%%%%%%%%%%%%%%%%%
%%%%%%%%%%%%%%%%%%%%%%%%%%%%%%%%%%%%%%%%%%%%%%%%%%%%%%%%%%%%%%%%%%%%%%%%%%%%%%%%%%%%%%%%%%%%%%%%%%%%%%%%%%%%
%%%%%%%%%%%%%%%%%%%%%%%%%%%%%%%%%%%%%%%%%%%%%%%%%%%%%%%%%%%%%%%%%%%%%%%%%%%%%%%%%%%%%%%%%%%%%%%%%%%%%%%%%%%%
%%%%%%%%%%%%%%%%%%%%%%%%%%%%%%%%%%%%%%%%%%%%%%%%%%%%%%%%%%%%%%%%%%%%%%%%%%%%%%%%%%%%%%%%%%%%%%%%%%%%%%%%%%%%
%%%%%%%%%%%%%%%%%%%%%%%%%%%%%%%%%%%%%%%%%%%%%%%%%%%%%%%%%%%%%%%%%%%%%%%%%%%%%%%%%%%%%%%%%%%%%%%%%%%%%%%%%%%%
\maketitle

\section{Introduction}
%\numsec

The main aim of this paper is to prove some results of the type of the Lyapunov Center Theorem for Hamiltonian systems. We consider the system 
\begin{equation}\label{eq:ham}
\dot{x}(t)=J_{2N}H'(x(t)),
\end{equation} 
 where $J_{2N}=\left[\ba{cc} 0 & Id_N \\ -Id_N & 0 \ea \right]$ is the standard symplectic matrix, $\Omega \subset \bR^{2N}$ is an open set containing the origin and $H \in C^2(\Omega, \bR)$ satisfy $H'(0)=0.$ 
Our goal is to study the existence of connected sets of closed trajectories of solutions  of the system \eqref{eq:ham}, emanating from the origin.
 
Obviously, this is a classical problem which goes back to the time of Poincar\'e, mainly with reference to the equations of celestial mechanics. The well-known result concerning the existence of periodic solutions is the celebrated Lyapunov Center Theorem proved in 1907, see \cite{Lyapunov}. It states, that if there exists a pair $\lambda_1 = \overline{\lambda_2}$ of purely imaginary eigenvalues of $J_{2N}H''(0)$, satisfying the so-called nonresonance condition, i.e. $\dfrac{\lambda_j}{\lambda_1} \notin \bZ$, where $\lambda_j, j=3, \ldots, 2N,$ are other eigenvalues of $J_{2N}H''(0),$ then there exists a connected branch of closed trajectories of solutions of the system, emanating from the origin.

 It is known, that if the nonresonance condition does not hold, it can happen that the system does not have nonstationary periodic solutions, see Example 9.2 in \cite{MawWill}. Therefore the problem of formulating conditions implying the existence of nonstationary periodic solutions in the neighbourhood of the origin has been studied by many researchers, among others Weinstein (see \cite{Weinstein}), Moser (see \cite{Moser}), Fadell and Rabinowitz (see \cite{FadRab}) and Bartsch (see \cite{Bartsch}). The assumptions of the theorems in these papers are given in terms of eigenvalues and eigenspaces of the matrix $J_{2N}H''(0)$. It is easy to observe that in all the above mentioned articles the considered eigenvalues have to be semisimple and the equilibrium is assumed to be nondegenerate.
The case of the isolated degenerate equilibrium has been considered by Szulkin, see \cite{Szulkin}, who has proved the existence of the sequence of solutions, converging to the equilibrium.

 The main result of our paper is Theorem \ref{thm:main}, which is of the type of the Lyapunov Center Theorem. We emphasize that it can be applied  also in the case of degenerate equilibrium without semisimple eigenvalues of $J_{2N}H''(0)$. Moreover, we obtain the existence of the connected set of solutions, not only the sequence, emanating from the origin.

The assumptions of our result are formulated with the use of the normal form of the matrix $J_{2N}H''(0)$, see condition \eqref{eq:sumaLap}. To formulate it we use the classification of normal forms of Hamiltonian matrices given in \cite{Chow} and count the numbers of indecomposable blocks of odd dimension, see Definition \ref{def:numbers}.  

We realize that verifying condition \eqref{eq:sumaLap} requires the determination of normal form of the matrix $J_{2N}H''(0)$, which in the general case can be complicated. However we observe that condition \eqref{eq:sumaLap} follows from  easier to verify assumptions of theorems of Lyapunov, Weinstein, Moser, Fadell and Rabinowitz, Bartsch, see assumptions (A0)-(A4) in Section \ref{ltt}. We emphasize that our theorem can be applied also in the situation when the assumptions of these theorems are not satisfied. Namely, in Example \ref{nonsemisimpleeven} we discuss the Hamiltonian system whose linearization does not have semisimple eigenvalues and assumptions of our theorem are fulfilled.

The second aim of our paper is to prove the global bifurcation theorem for parametrized Hamiltonian systems, see Theorem \ref{gbthoe}. We claim that the condition similar to \eqref{eq:sumaLap} implies the global bifurcation of nonstationary $2\pi$-periodic solutions of such system.

Our article is organised in the following way. In Section \ref{sec:preliminaries} we introduce the necessary notation and definitions. The main results of our paper, i.e. Theorems \ref{thm:main} and \ref{gbthoe}, are formulated in Section \ref{sec:main}. Moreover, we discuss there the assumptions of these theorems. Their proofs are given in Section \ref{sec:proof}.

%The problem of the existence of periodic solutions in the neighbourhood of the equilibrium is related to the problem of bifurcation of solutions of the parameterized system. For such system the condition similar to \eqref{eq:sumaLap} implies the global bifurcation of nonstationary $2\pi$-periodic solutions. The second main result of our paper is the global bifurcation theorem, see Theorem \ref{gbthoe}.

%It is known, that if the nonresonance condition does not hold, then it can happen, that the system does not have non-constant periodic solutions, see \cite{MawWill}. Therefore, it is important to formulate conditions for continuation of periodic solutions, weaker than nonresonance condition. \textcolor{blue}{opisac twierdzenia Weinsteina, Mosera, Fadella-Rabinowitza, Szulkina, Bartscha}

%The assumptions of all the theorems mentioned above are given in terms of eigenvalues and eigenspaces of matrix $J_{2N}H''(0)$. However, it is easy to observe that in all cases the considered eigenvalues have to be semisimple.

%The main aim of our paper is to prove the result generalizing above theorems to the case of nonsemisimple as well as semisimple purely imaginary eigenvalues. 

%--------------------------------------------------------------
%--------------------------------------------------------------

\section{Basic definitions}\label{sec:preliminaries}

The aim of this section is to introduce the notation that we use in our article. In particular, we discuss the concept of the normal form of Hamiltonian matrix and   describe the notion of local and global bifurcation of nonstationary periodic solutions of autonomous Hamiltonian systems.

\subsection{Normal forms of Hamiltonian matrices}\label{subs:normal}
One of the main tools that we use to prove our results is the theory of real normal forms of Hamiltonian matrices.  From now on all the considered matrices are assumed to be real.

 The theory of normal forms of linear Hamiltonian systems has been developed by many authors, see discussion in the Bibliographical Notes in Chapter 2 of \cite{Chow}. In our paper we use the approach presented in Chapter 2.6 of \cite{Chow}. Below we recall some elements of this theory.

\bdf
A $(2N \times 2N)$-matrix $S$ is called symplectic, if $S^TJ_{2N}=J_{2N}S^{-1}$ or equivalently $S^{-1}J_{2N}=J_{2N}S^{T}.$
A $(2N \times 2N)$-matrix $M$ is called Hamiltonian (or infinitesimally symplectic) if $M^{T}=J_{2N}MJ_{2N}.$
\edf

\br It is known that the matrix $M$ is Hamiltonian iff $M=J_{2N} A$, where $A$ is a symmetric matrix.
\er

\bdf
The Hamiltonian matrix $M$ defined on a symplectic vector space $(V,J_{2n})$ is called decomposable if $V=V_1\oplus V_2$, where $V_1$ and $V_2$ are proper, $M$-invariant and $J_{2N}$-orthogonal  (i.e. $v_1^T J_{2N}v_2=0$ for all $v_1 \in V_1, v_2 \in V_2$) symplectic subspaces of $V$. $M$ is called indecomposable if $A$ is not decomposable.
\edf

\bdf
Two  Hamiltonian $(2N \times 2N)$-matrices $M_1$ and $M_2$ are called symplectically similar if there exists a symplectic matrix $S$ such that $M_2=S^{-1}M_1S.$
\edf

%From now on we assume that $M$ is a fixed Hamiltonian matrix.

\begin{Remark} \label{rem:normalform} Obviously, the symplectic similarity is an equivalence relation. Moreover, using Theorem 6.13 of \cite{Chow} we obtain that each Hamiltonian $(2N \times 2N)$-matrix $M$ is symplectically similar to the matrix of the form 
\begin{equation}\label{eq:Mj}
\left[ \ba{ll} \diag(D_1,\ldots, D_s) & \diag (B_1,\ldots, B_s) \\ \diag (C_1,\ldots,C_s )& \diag (-D_1^T,\ldots,-D_s^T) \ea  \right],
\end{equation}
where 
\begin{enumerate}[(1)]
\item $B_j,C_j,D_j$ are $(N_j \times N_j)$-matrices, $B_j^T=B_j, C_j^T=C_j, j = 1,\ldots,s$,
\item $N_1 + \ldots + N_s = N$,
\item $M_j=\left[ \begin{array}{rr} D_j &B_j\\ C_j &-D_j^{T}\end{array}\right]$ is one of the Hamiltonian matrices listed on pages 118-122 of \cite{Chow}. In particular, $M_j$ is indecomposable. %in the sense of Definition 6.12 of \cite{Chow}. 
\end{enumerate}
The matrix in the form \eqref{eq:Mj}, being in the equivalence class of matrix $M$, is called the normal form of $M$. 
\end{Remark}

The condition that $M_j$ is indecomposable Hamiltonian matrix implies that its spectrum can have only  one of the following forms: $\{0\}, \{\pm \alpha\}, \{\pm \beta i\}, \{\pm \alpha \pm \beta i\}$ and the elements of such spectrum are the eigenvalues of  $M$. The most important situation from the point of view of our computations is the case  of matrices $M_j$ such that $\sigma(M_j)=\{\pm \beta_0 i\}, \beta_0 > 0$. For such matrices the normal forms depend on the parity of $N_j$. Moreover, both in the cases of $N_j$ odd and even, there are two possible normal forms. Following the lists of such forms given in \cite{Chow}, we consider four  possible forms of the matrix $M_j$ such that $\sigma(M_j)=\{\pm \beta_0 i; \beta_0>0\}$, which can be described in the two following cases:

\be[(C1)]
\item $N_j$ is odd. In this case there are two possible normal forms, for $\epsilon = \pm 1$,  given by $M_{j,\epsilon}=J_{2N_j}H_{\epsilon}''(0)$ where
\begin{equation*}
\begin{split}
&H_{\epsilon}(x_1,\ldots,x_{N_j},y_1,\ldots,y_{N_j})=\\&= -\epsilon \beta_0 \sum_{i=1}^{\frac{N_j-1}{2}} (-1)^{i+1} (x_i x_{N_j+1-i} + y_i y_{N_j+1-i}) + \sum_{i=1}^{N_j-1} x_i y_{i+1} + \frac{1}{2} \epsilon \beta_0(-1)^{\left[\frac{N_j}{2} \right] +1} \left(x^2_{\frac{N_j+1}{2}} + y^2_{\frac{N_j+1}{2}}\right).
\end{split}
\end{equation*}

Then
\beq \label{c2} M_{j,\epsilon}=\left[\ba{ccccc;{2pt/2pt}ccccc}

0 &  & &  &   &  &   & &  & -\epsilon \beta_0 

\\ 1 & 0 &  &  &  &  &  & &\epsilon \beta_0 &

\\ &  \ddots &0 & &  &  & &  \udots &  &  

\\ &  &   & \ddots  & & & \epsilon \beta_0 &  &  &  

\\   & &  & 1 & 0 & - \epsilon \beta_0 &  &   &  &  
\\ \hdashline
  &  &  &  & \epsilon \beta_0 & 0 & -1 &  & &

\\  &  & & - \epsilon \beta_0 & &  & 0&-1 &  & 

\\  & &  \udots &  &  & &  &  &\ddots & 

\\    & - \epsilon \beta_0 &  &  & &  &  & & \ddots &-1 

\\    \epsilon \beta_0 &  &  & &  &  & &  & & 0
\ea \right].
\eeq
For $N_j=1$ this matrix coincide with the one given in List I (3) on page 119 of \cite{Chow}, while in the case $N_j>1$ with the one in List II (5) on page 121.

%$M_{j}$ is semisimple with $\sigma(M_{j})=\{\pm  \beta_0 i; \beta_0>0\}$. In this case there are two possible normal forms, given by $M_{j, \epsilon}=J_{2}H_{\epsilon}''(0)$ for $H_{\epsilon}(x,y)=\frac{1}{2} \epsilon \beta_0(x^2+y^2),$ for $\epsilon = \pm 1.$ 

%Then 
%\beq \label{c1} M_{j,\epsilon}=\left[\ba{cc}
%0 &\epsilon \beta_0\\
%-\epsilon \beta_0 &0
%\ea\right],\eeq
% see List I (3) on page 119 of \cite{Chow}.

%\item $M_{j}$ is nonsemisimple, $N_j$ is odd and $\sigma(M_{j})=\{\pm  \beta_0 i; \beta_0>0\}$. In this case 

\item $N_j$ is even. In this case there are two possible normal forms, for $\epsilon = \pm 1$,  given by $M_{j,\epsilon}=J_{2N_j}H_{\epsilon}''(0)$ where

\begin{equation*}
H_{\epsilon}(x_1,\ldots,x_{N_j},y_1,\ldots,y_{N_j})= \beta_0 \sum_{i=1}^{N_j/2} (x_{2i-1}y_{2i} - x_{2i}y_{2i-1}) + \sum_{i=1}^{N_j-2} x_i y_{i+2} - \frac{1}{2} \epsilon (x_{N_j-1}^2 + x_{N_j}^2).
\end{equation*}

Then
\beq \label{c3} M_{j,\epsilon}=\left[\ba{cccc;{2pt/2pt}cccc}

A_2 &  &  &  &   &  &   & 

\\ I_2 & \ddots &  &  &  &  &  &  

\\  & \ddots & \ddots&   & &  & &     

%\\ \vdots& \ddots &  \ddots & \ddots  & \vdots & 0 & \epsilon \beta_0 &\udots  &  & \vdots 

\\    &  & I_2 & A_2 &  &  & &   
\\ \hdashline
  &  &  &   & A_2 & -I_2 &  & 

\\  &  & &  & & \ddots & \ddots  & 

%\\ \vdots & &  \udots &  &  & &  & 0 &\ddots & 0

\\   &  &  & &  &  &\ddots & -I_2 

\\     &  & & \epsilon I_2 &  &   & & A_2
\ea \right],
\eeq
where $I_2=\left[\ba{cc}
1 & 0\\
0  & 1
\ea\right]$, $A_2=\left[\ba{cc}
0 &-\beta_0\\
\beta_0 & 0 
\ea\right],$  see List II (4) on page 120 of \cite{Chow}.
\ee

\begin{Remark}
Note that for a given indecomposable Hamiltonian matrix with spectrum $\{\pm \beta_0 i\}$ the normal form can be explicitely determined with the use of the construction given in \cite{Burgoyne}.
\end{Remark}

\begin{Remark}
We are aware that one can consider also another normal forms for matrices with purely imaginary eigenvalues, see for example \cite{LaubMeyer, MeyerHallOffin}. Obviously, for a given Hamiltonian matrix, its normal forms presented in \cite{Chow} and \cite{LaubMeyer, MeyerHallOffin} are symplectially similar.
\end{Remark}

Below we classify purely imaginary eigenvalues of a Hamiltonian matrix $M$.

 \begin{Definition}
The eigenvalue $\pm \beta_0 i \in \sigma (M)$ is called
\begin{enumerate}
\item[1)] simple if there is only one block $M_j$ corresponding to $\pm \beta_0 i$ in the decomposition \eqref{eq:Mj} and $M_j$ is of dimension two,  

\item[2)] semisimple if all the blocks $M_j$ corresponding to $\pm \beta_0 i$ in the decomposition \eqref{eq:Mj} are of dimension two,

\item[3)] partially semisimple if in the decomposition \eqref{eq:Mj} there are at least two blocks $M_{j_1},M_{j_2}$ corresponding to $\pm \beta_0 i$   of dimension two and greater than two, respectively,

\item[4)] strictly nonsemisimple  if the dimension of any block $M_j$ in the decomposition \eqref{eq:Mj} corresponding to $\pm \beta_0 i$ is greater than  two.
\end{enumerate}
\end{Definition}

\br Let $\cC(M)$ be the complex Jordan normal form of the matrix $M$ and let  $\beta_0 i \in  \sigma(M)$. Then
\be
\item the eigenvalue $\beta_0 i$ is simple iff there is exactly one elementary Jordan block of the matrix $\cC(M)$ corresponding to $\beta_0 i$ and it is of dimension one i.e. the algebraic multiplicity  of  $\beta_0 i$ equals one,

\item  the eigenvalue $\beta_0 i$ is semisimple iff the dimension of any elementary Jordan block of the matrix $\cC(M)$ corresponding to $\beta_0 i$ is equal to one i.e.     the algebraic multiplicity of $ \beta i$ equals its geometric multiplicity,

\item  the eigenvalue $\beta_0 i$ is partially semisimple iff there are elementary Jordan blocks of the matrix $\cC(M)$ corresponding to $\beta_0 i$ of dimension one as well as greater than one i.e. the algebraic multiplicity of $\beta_0 i$ is different from its geometric multiplicity and there is at least one elementary Jordan block of the matrix $\cC(M)$  corresponding to $\beta_0 i$ of dimension one,

\item  the eigenvalue $\beta_0 i$ is strictly nonsemisimple iff the dimension of any elementary block of the matrix $\cC(M)$ corresponding to $\beta_0 i$ is greater than one.
\ee
\er
Our main results (see Section \ref{sec:main}) are formulated with the use of the number of blocks described in case (C1), i.e. the number of the blocks with $N_j$ odd. Additionally, we distinguish situations of even and odd value of $\dfrac{N_j+1}{2}$, and moreover of positive and negative value of $\epsilon.$ Below we introduce   four numbers depending on these values.

%Our main result depends on the numbers of the blocks described above, with $N_j$ odd. However, in section \ref{sec:proof} we observe that the cases of $\dfrac{N_j+1}{2}$ odd and even should be treated separately. Moreover, we distinguish cases of positive and negative value of $\epsilon.$ Below we introduce the notation for  these numbers. 

\begin{Definition} \label{def:numbers}
Fix $\pm \beta_0 i \in \sigma(M)$. Denote   by
\begin{enumerate}
\item $o_{\pm}(\beta_0,M)$  the number of $(2N_j \times 2N_j)$-blocks $M_{j,\epsilon}$  in the decomposition \eqref{eq:Mj}  such that   $N_j$ is odd, $\frac{N_j+1}{2}$ is odd and   $\epsilon=\pm 1$,
\item $e_{\pm}(\beta_0,M)$  the number of $(2N_j \times 2N_j)$-blocks $M_{j,\epsilon}$  in the decomposition \eqref{eq:Mj}  such that $N_j$ is odd, $\frac{N_j+1}{2}$ is even and  $\epsilon=\pm 1$.
\end{enumerate}

\end{Definition}

%\begin{Definition}\label{rem:numbers} Fix $\pm \beta_0 i \in \sigma(J_{2N}A)$. Denote:
%\begin{enumerate}
%\item by $o_+(\beta_0,J_{2N}A)$  the number of $(2N_j \times 2N_j)$-blocks $M_j$ (corresponding to the eigenvalue $\pm \beta_0 i$)   such that $N_j$ is odd and $\frac{N_j+1}{2}$ is odd and the value of $\epsilon$ is equal to $+1$,
%\item by $o_-(\beta_0,J_{2N}A)$  the number of $(2N_j \times 2N_j)$-blocks $M_j$ (corresponding to the eigenvalue $\pm \beta_0 i$)   such that $N_j$ is odd and $\frac{N_j+1}{2}$ is odd and the value of $\epsilon$ is equal to $-1$,
%\item by $e_+(\beta_0,J_{2N}A)$  the number of $(2N_j \times 2N_j)$-blocks $M_j$ (corresponding to the eigenvalue $\pm \beta_0 i$)   such that $N_j$ is odd and $\frac{N_j+1}{2}$ is even and the value of $\epsilon$ is equal to $+1$,
%\item by $e_-(\beta_0,J_{2N}A)$  the number of $(2N_j \times 2N_j)$-blocks $M_j$ (corresponding to the eigenvalue $\pm \beta_0 i$)   such that $N_j$ is odd and $\frac{N_j+1}{2}$ is even and the value of $\epsilon$ is equal to $-1$,

%\end{enumerate}
%\end{Definition}

 \br \label{rem:semi} Note that if $M_{j,\epsilon}$ is one of the blocks  corresponding to the semisimple eigenvalue $\pm \beta_0 i$ (and therefore $N_j=1$), then $\frac{N_j+1}{2}=1$ is odd. Therefore  the number of blocks corresponding to the semisimple eigenvalue $\pm \beta_0 i$ equals $o_+(\beta_0, M) +o_-(\beta_0, M)$ and $e_+(\beta_0,M)=e_-(\beta_0,M)=0$. However, we emphasize that we do not restrict our attention only to the symplectically diagonalizable case i.e. when all the blocks corresponding to  purely imaginary eigenvalues are of dimension two.
\er

\subsection{Global bifurcation} The method that we use to study the existence of periodic solutions of the system \eqref{eq:ham} is to change the problem to bifurcation one and study the phenomenon of global bifurcation. In this subsection we introduce the notion of bifurcation. 

Consider a Hamiltonian $H \in C^2(\Omega, \bR)$, where $\Omega \subset \bR^{2N}$ is an open subset, satisfying condition $(H')^{-1}(0) \cap \Omega =\{s_1, \ldots, s_k\}.$

We investigate the existence of $2\pi$-periodic solutions of the following family of Hamiltonian systems:

\begin{equation}\label{eq:parham}
\dot{x}(t)=\lambda J_{2N}H'(x(t)),
\end{equation}
where $\lambda>0.$

Put $C_{2\pi}(\Omega)=\{u \in C([0,2\pi], \Omega); u(0)=u(2\pi)\}.$ We consider $2\pi$-periodic solutions of the system \eqref{eq:parham} as elements  $(u(t), \lambda) \in \Comega\times (0, +\infty).$ Moreover, identifying $s_1, \ldots, s_k \in \Omega$ with constant functions in $\Comega$, we define two subsets $\cT, \cN \subset \Comega \times (0,+\infty)$ as follows 
\begin{equation*}
\begin{split}
&\cT=\{s_1, \ldots, s_k\}\times(0, +\infty),\\
&\cN=\{(u(t),\lambda): (u(t),\lambda) \text{ is a nonstationary } 2\pi\text{-periodic solution of \eqref{eq:parham}}\}.
\end{split}
\end{equation*}

The elements of $\cT$ are called trivial solutions of the system \eqref{eq:parham} and elements of $\cN$  nontrivial ones.

Fix $(s_{k_0}, \lambda_0) \in \cT$. If $(s_{k_0}, \lambda_0) \in \cT \cap \cl(\cN)$  denote by $\cC(s_{k_0}, \lambda_0)$ the continuum (i.e. closed connected component) of $\cl(\cN)$ containing $(s_{k_0}, \lambda_0)$, where the closure of $\cN$ is taken in $C_{2\pi}(\Omega) \times [0, +\infty).$ For $(s_{k_0}, \lambda_0) \notin \cT \cap \cl(\cN)$ we define $\cC(s_{k_0}, \lambda_0)=\emptyset.$

\begin{Definition}
A point $(s_{k_0}, \lambda_0)\in \cT$ is called
\be 
\item a local bifurcation point of nonstationary $2\pi$-periodic solutions of the system \eqref{eq:parham} if $(s_{k_0}, \lambda_0) \in \cl(\cN)$, i.e. $\cC(s_{k_0}, \lambda_0) \not =\emptyset$,
\item a branching point  of nonstationary $2\pi$-periodic solutions of the system \eqref{eq:parham} if $\cC(s_{k_0}, \lambda_0) \not = \{(s_{k_0}, \lambda_0)\},$

\item a global bifurcation point of nonstationary $2\pi$-periodic solutions of the system \eqref{eq:parham} if either the continuum $\cC(s_{k_0}, \lambda_0)$ is not compact in $\Comega\times(0, +\infty)$ or it is compact and $\cC(s_{k_0}, \lambda_0) \cap (\cT \setminus \{(s_{k_0}, \lambda_0)\}) \neq \emptyset.$
\ee
\end{Definition}

\begin{Remark}
In other words, by a local bifurcation phenomenon  we understand the existence of a sequence $\{(x_p(t), \lambda_p)\} \subset \cN$ converging (in $L_{\infty}$- norm) to $(s_{k_0}, \lambda_0)$, whereas by a branching phenomenon we understand the existence of a closed connected set of nontrivial solutions containing $(s_{k_0}, \lambda_0)$. On the other hand, when a global bifurcation phenomenon occurs, this emanating set either is not compact or it is compact and  touches the set of stationary solutions in more than one point. Obviously, a global bifurcation point is a branching point. Moreover, it is evident that  a branching point is a local bifurcation point.
\end{Remark}

The necessary condition for the existence of a  local bifurcation  of nonstationary $2\pi$-periodic solutions of the system \eqref{eq:parham} is well-known, see for example  Remark 3.4 of \cite{DanRyb}. For the completeness of the paper, we recall it below. Put
$$\Lambda(s_{k_0}) =\left\{\frac{m}{\beta} \in (0, +\infty): \pm i\beta \in \sigma(J_{2N}H''(s_{k_0})), \beta>0, m \in \bN \right\}.$$

\begin{Fact}
If $(s_{k_0}, \lambda_0)\in \cT$ is a local bifurcation point of nonstationary $2\pi$-periodic solutions of the system  \eqref{eq:parham}, then $\lambda_0 \in \Lambda(s_{k_0}).$
\end{Fact}

To study the sufficient condition of global bifurcation we use the theorem proved by Dancer and the second author in \cite{DanRyb}. This result is formulated with the use of the bifurcation index. Below we recall this notion.

Let $j \in \bN, s_i \in \{s_1, \ldots, s_k\}$ and $\lambda \in \bR$. Define a $(4N \times 4N)$-matrix $T_j^N(\lambda,H''(s_i))$ by 
\begin{equation}\label{eq:Tj}
T_j^N(\lambda H''(s_i))=\left[ \begin{array}{cc} -\frac{\lambda}{j}H''(s_i) & J_{2N}\\
-J_{2N} &-\frac{\lambda}{j}H''(s_i)
\end{array}\right]
\end{equation}

Choose $\mu>0$ such that $[\lambda_0 - \mu, \lambda_0 + \mu] \cap \Lambda(s_{k_0}) =\{\lambda_0\}.$ Note that since the set $\Lambda(s_{k_0})$ does not have finite accumulation points, such a choice is possible. 

Define a bifurcation index 
$\ds  \bif(s_{k_0}, \lambda_0) \in \bigoplus_{j=1}^{\infty} \bZ $ as follows  

\beq \label{bifin} \bif(s_{k_0}, \lambda_0)=(\eta_1(s_{k_0}, \lambda_0), \ldots, \eta_j(s_{k_0}, \lambda_0), \ldots) \in  \bigoplus_{j=1}^{\infty} \bZ ,\eeq
where 
\beq \label{bifi} \eta_j(s_{k_0}, \lambda_0)=\ib(s_{k_0}, H')\cdot (m^-(T_j^N((\lambda_0+\mu) H''(s_{k_0})))-m^-(T_j^N((\lambda_0-\mu) H''(s_{k_0})))).\eeq
In the above formula $\ib(s_{k_0}, H')$ stands for the Brouwer index of $s_{k_0}$  and $\morse(A)$ denotes the Morse index of a symmetric matrix $A$. 

Moreover, set  $\ds \Theta = (0, \ldots, 0, \ldots) \in \bigoplus_{j=1}^{\infty} \bZ.$

In the theorem below we describe properties of continua of nonstationary periodic solutions of autonomous Hamiltonian systems. This theorem  is a reformulation of Theorem 3.3 of \cite{DanRyb}. In particular, to make the formulation of the theorem more readable, we use the space of continuous functions rather than the Sobolev space $H^1_{2\pi}$.

\begin{Theorem} \label{gbth} Fix $(s_{k_0},\lambda_0) \in \{s_{k_0}\} \times \Lambda(s_{k_0})$ satisfying  $\bif(s_{k_0}, \lambda_0) \neq \Theta.$ Then $(s_{k_0}, \lambda_0) \in \cT$ is a global bifurcation point of nonstationary $2\pi$-periodic solutions of the system \eqref{eq:parham}. Moreover, if the continuum $\cC(s_{k_0}, \lambda_0)$ is compact in $\Comega \times (0, +\infty)$, then
\begin{enumerate}
\item[(a)] $\cC(s_{k_0}, \lambda_0) \cap \cT$ is finite,
\item[(b)] $\displaystyle{\sum_{(\hat{s}, \hat{\lambda})\in \cC(s_{k_0}, \lambda_0) \cap \cT}} \bif(\hat{s},\hat{\lambda})=\Theta.$
 
\end{enumerate}
\end{Theorem}

\br \label{drv}
To prove Theorem 3.3 of \cite{DanRyb} we have considered the $2 \pi$-periodic solutions of the family  \eqref{eq:parham} as critical points of a family of $\sone$-invariant $C^2$-functionals $\Phi : H^1_{2 \pi}\times (0,+\infty) \to \bR$ defined by the following formula
$$
\Phi(x,\lambda)=\int_0^{2\pi} (-J \dot x(t),x(t))-\lambda H(x(t)) dt .
$$
In other words,   to study properties of sets of $2 \pi$-periodic solutions of the family  \eqref{eq:parham} it is enough to study the set of solutions of the following equation  $$\nabla_x\Phi(x,\lambda)=0.$$

Applying the Amann-Zehnder reduction, see \cite{[AMZE]}, we have reduced the problem to a finite-dimensional one with $\sone$-symmetries and variational structure. Finally we have applied the finite-dimensional version of the degree for $\sone$-equivariant gradient maps, see \cite{Geba}, to prove the global bifurcation theorem for $\sone$-equivariant gradient maps.

In 2011, see \cite{[GORY]}, we defined the infinite-dimensional generalization of  the degree for $\sone$-equivariant gradient maps. That is why we can prove Theorem 3.3 of \cite{DanRyb} without the Amann-Zehnder reduction, working directly in infinite-dimensional spaces. The proof would be in spirit the same as the proof of the famous theorem due to Rabinowitz, see \cite{Rabinowitz}, where the Leray-Schauder degree was applied as the topological tool to prove the existence of phenomenon of global bifurcation.

\er

\br \label{kb}
It is worth to point out that, in 1953, Krasnosel'skii used variational methods to prove that if a Fr\'echet differentiable
operator $A : H \to H$ is the gradient of a weakly continuous functional in a Hilbert space $H$
and A(0) = 0, then each characteristic value of
$A'(0)$ is a local bifurcation point of the equation $A(x)=\lambda x,$ see \cite{CH}, \cite{Krasnosielski}, \cite{KZ}. Using the Conley index, Ize also proved this theorem, see \cite{Ize}.
This deep result has stimulated a lot of contributions in variational bifurcation theory. The following question is important: Is there a continuum  of nontrivial solutions emanating from each characteristic value of $A'(0)$? The answer is No. Relevant examples were discussed in \cite{Ambrosetti}, \cite{Bohme}, \cite{Ize}, \cite{Marino}, \cite{Takens}.

\er

\br Note that 
\be

\item Although the situation in Theorem 3.3 of \cite{DanRyb} is slightly different, we claim that Theorem \ref{gbth} is a reformulation of this result. Let us first consider the case $\Omega = \bR^{2N}$. In Theorem 3.3 of \cite{DanRyb}, a global bifurcation is considered in the space $H^1_{2\pi} \times (0, +\infty)$, where the Hilbert space $H^1_{2\pi}$ is endowed with a scalar product $\langle u, v \rangle_{H^1_{2\pi}} = \displaystyle{ \int_0^{2\pi} (\dot{u}(t), \dot{v}(t))+(u(t), v(t)) dt}$. However, since we study the continua of solutions of the system \eqref{eq:parham}, we can formulate this result in the space of continuous functions $C_{2\pi}(\bR^{2N})$ with the supremum norm $\|\cdot \|_{\infty}$. First of all,  we consider the set $\cN\subset H^1_{2\pi} \times(0, +\infty) \subset C_{2\pi}(\bR^{2N})\times(0,+\infty)$ of nontrivial solutions ot the system \eqref{eq:parham}. It is worth to point out, that the closure of $\cN$ in both spaces consists of solutions of \eqref{eq:parham}, and consequently these closures coincide.

Moreover, the compactness of the continuum $\cC(s_{k_0}, \lambda_0)$ in the space $H^1_{2\pi}\times (0, +\infty)$ is equivalent to the compactness in $C_{2\pi}(\bR^{2N})\times(0,+\infty)$. Indeed, since the embedding $H^1_{2\pi} \subset C_{2\pi}(\bR^{2N})$ is continuous (see Proposition 1.1. of \cite{MawWill}), if the set $\cC(s_{k_0}, \lambda_0)$ is compact in $H^1_{2\pi}\times (0,+\infty)$, then it is also compact in $C_{2\pi}(\bR^{2N})\times (0, +\infty)$. 

To prove the latter implication, we take the sequence $\{(u_n, \lambda_n) \} \subset \cC(s_{k_0}, \lambda_0)\subset \cl(\cN) \subset H^1_{2\pi} \times (0, +\infty) \subset C_{2\pi}(\bR^{2N}) \times(0,+\infty).$ If $\cC(s_{k_0}, \lambda_0)$ is compact in $C_{2\pi}(\bR^{2N}) \times (0, +\infty)$, then we can choose the subsequence $\{(u_{n_k}, \lambda_{n_k})\}$ convergent to $(u_{\star}, \lambda_{\star}) \in C_{2\pi}(\bR^{2N}) \times (0, +\infty).$ We will show that $\|(u_{n_k}, \lambda_{n_k})-(u_{\star}, \lambda_{\star})\|_{H^1_{2\pi}\times (0, +\infty)} \to 0.$ This is a consequence of the fact that continuum $\cC(s_{k_0}, \lambda_0)$ consists of solutions, and $H'$ is a $C^1$-function. Therefore,  $$\|\dot{u}_{n_k}-\dot{u}_{\star}\|_{\infty}=\|\lambda_{n_k} J_{2N}H'(u_{n_k}(\cdot))-\lambda_{\star}J_{2N}H'(u_{\star}(\cdot))\|_{\infty} \to 0$$ and in consequence, since $\|u\|_{L^2} \leq 2\pi \|u\|_{\infty},$ 
$$\|u_{n_k}-u_{\star}\|^2_{H^1_{2\pi}} = \|\dot{u}_{n_k} -\dot{u}_{\star}\|^2_{L^2}+\|u_{n_k}-u_{\star}\|^2_{L^2}\to 0,$$ hence $\|(u_{n_k}, \lambda_{n_k})-(u_{\star}, \lambda_{\star})\|_{H^1_{2\pi}\times (0, +\infty)} \to 0.$

\item Since the celestial mechanics provides a variety of examples, where the Hamiltonian $H$ is not defined on the whole space $\bR^{2N}$, we consider the space $\Comega$ instead of $C_{2\pi}(\bR^{2N})$. It is easy to observe, that our assertion remains valid in this case,

\item  the continuum $\cC(s_{k_0},\lambda_0)$ does not reach the level $\lambda=0,$ i.e. $\cC(s_{k_0},\lambda_0) \cap (\Comega \times \{0\}) = \emptyset$, see  Step 2 of the proof of Theorem  3.3 in \cite{DanRyb},  

\item if $s_{k_0}$ is a nondegenerate critical point of $H$, then computing the bifurcation index reduces to a linear algebra problem; otherwise we additionally have to compute the Brouwer index $\ib(s_{k_0}, H')$ of a degenerate critical point $s_{k_0}$, which is not an easy  task, 

\item condition $(b)$ of the above theorem implies that in some cases one can prove that the continuum $\cC(s_{k_0},\lambda_0)$ is not compact by showing that all the nontrivial coordinates of bifurcation indices have the same sign.
\ee
\er

\section{Main results}\label{sec:main}
In this section we formulate and discuss the main results of our paper. The proofs of these theorems are postponed to the next section.  

In particular, in Subsection \ref{ltt} we formulate  a generalization of the Lyapunov Center Theorem, see Theorem \ref{thm:main}, and present an example, see Example \ref{nonsemisimpleeven}, of Hamiltonian system with strictly nonsemisimple eigenvalue of even multiplicity, for which assumptions of this theorem are fulfilled. Moreover, we show that assumptions of some famous   theorems imply the ones  of Theorem \ref{thm:main}, see Propositions \ref{colct} -  \ref{cofrmb} and Remarks \ref{general}, \ref{szulk}. Additionally, we show that these classical theorems  cannot  be applied  to the study of periodic solutions of the Hamiltonian system considered in Example \ref{nonsemisimpleeven}.
In Subsection \ref{secgbth} we  formulate a global bifurcation theorem for autonomous Hamiltonian systems and discuss its assumptions.

\subsection{Lyapunov-type theorem} \label{ltt}

In this section we formulate results concerning the existence of nonstationary periodic solutions in a neighborhood of a stationary one of the following Hamiltonian system 

\beq \label{has}
\dot x(t) =J_{2N} H'(x(t)),
\eeq
 where $\Omega \subset \bR^{2N}$ is open, $H \in C^2(\Omega,\bR)$ and $0 \in \Omega$ is an isolated critical point of the Hamiltonian $H$.

The following theorem is a generalization of the Lyapunov Center Theorem. The nonresonance condition, see condition $\mathbf{(A0)}$ in the discussion below, is replaced by the more general condition \eqref{eq:sumaLap}. Recall, that the numbers $o_{\pm}, e_{\pm}$ have been introduced in the Definition \ref{def:numbers} and $\ib$ denotes the Brouwer index.

\begin{Theorem} \label{thm:main}
For an open subset  $\Omega \subset \bR^{2N}$ 
consider a Hamiltonian  $H \in C^2(\Omega, \bR)$ such that $0 \in \Omega$ is an isolated critical point of $H$ and  $\ib(0,H') \neq 0$. If there exists a pair $\pm  \beta_0 i$ of purely imaginary eigenvalues of the matrix  $J_{2N}H''(0)$ such that  
\begin{equation} \label{eq:sumaLap}
o_+(\beta_0,J_{2N}H''(0))-o_-(\beta_0,J_{2N}H''(0))-e_+(\beta_0,J_{2N}H''(0))+e_-(\beta_0,J_{2N}H''(0)) \neq 0,
\end{equation}
  then there exists a connected set  $\cC_{\beta_0} \subset \Omega$ of closed trajectories of solutions of the system \eqref{has} emanating from the origin. Moreover
\beq \label{crezo}  \forall  \epsilon > 0 \: \exists \delta > 0 \: \forall u_0 \in \cC_{\beta_0}\: \text{ if }  \|u(\cdot, u_0)\|_{L^\infty} < \delta \text{ then } \mid T(u(\cdot, u_0)) - \frac{2 \pi}{\beta_0}\mid < \epsilon,\eeq where $u(\cdot, u_0)$ is a solution of the system \eqref{has} with the initial condition $x(0)=u_0$ and  $T(u(\cdot, u_0))$ is a period  of $u(\cdot, u_0)$.
\end{Theorem}

 \br Note that the period $T(u(\cdot, u_0))$   in the above theorem doesn't have to be minimal. However, we can obtain the result concerning minimal periods of solutions assuming an additional nonresonance condition. Let $\cB=\{\beta_1, \beta_2, \ldots, \beta_k\}$ be such that $\beta_1>\beta_2>\ldots>\beta_k>0$
and $\sigma(J_{2N} H''(0)) \cap i \bR =\{\pm \beta_1 i, \ldots, \pm \beta_k i\}$. We say that a nonresonance condition for $\beta_{k_0} \in \cB$ holds if $\frac{\beta_{j}}{\beta_{k_0}} \notin \bZ$ for $j=1, \ldots, k_0-1$. Assuming additionally in Theorem \ref{thm:main} that the nonresonance condition holds for $\beta_0$ we obtain  the condition \eqref{crezo}  with  $T(u(\cdot, u_0))$  being the minimal period of $u(\cdot, u_0)$. 
\er

\begin{Corollary}
The above remark allows us to distinguish, in a small neighbourhood of the origin, connected sets of solutions by the  minimal periods. Using this fact we can obtain the following result concerning the multiplicity of sets of solutions. Let $\cB$ be the set  defined as  in the above remark. Denote by $\cB_n$ the set of numbers $\beta_l \in \cB$ such that the nonresonance condition for $\beta_l$ holds. Then in a sufficiently small neighbourhood of $0 \in \bR^{2N}$ there exist at least $\card \cB_n$ of connected sets of non-constant periodic orbits of the system \eqref{has} emanating from the origin.
\end{Corollary}

 \br Note that the condition \eqref{eq:sumaLap} can be verified directly if we know the normal form of the matrix $J_{2N} H''(0)$. In the example below we illustrate Theorem \ref{thm:main}. Notice that in this example the spectrum of  $J_{2N}H''(0)$ is  the set $\{\pm \beta_0 i\}$ and the eigenvalue $\beta_0 i$ is strictly nonsemisimple and of even multiplicity. Obviously, a similar example with partially semisimple eigenvalue can be easily constructed.
\er

\bex \label{nonsemisimpleeven}
Let   $H \in C^2(\bR^{2N},\bR)$ be such that $H(0)=0$ and $H'(0)=0$, i.e. $H(x)=\frac{1}{2}\langle H''(0)x, x\rangle+h.o.t.$ Assume moreover that $\sigma(J_{2N} H''(0))=\{\pm \beta_0 i\},$ where $\beta_0>0$, and  $J_{2N} H''(0)$ is in the normal form \eqref{eq:Mj} with $s=3, N_1=4k+1,N_2=4k-1, k \in \bN,$ $N_3$ even and  $N=N_1+N_2+N_3$.  
Therefore,   $J_{2N} H''(0)$   is of the form:

$$J_{2N} H''(0)=\left[ \ba{ll} \diag(D_{1,\epsilon}, D_{2,\epsilon},D_{3,\epsilon}) & \diag (B_{1,\epsilon},B_{2,\epsilon},B_{3,\epsilon}) \\ \diag (C_{1,\epsilon}, C_{2,\epsilon},C_{3,\epsilon} )& \diag (-D_{1,\epsilon}^T,-D_{2,\epsilon}^T,-D_{3,\epsilon}^T) \ea  \right],$$

\be

\item $B_{j,\epsilon}, C_{j,\epsilon}, D_{j,\epsilon}$ are $(N_j \times N_j)$-matrices, $B_{j,\epsilon}^T=B_{j,\epsilon}, C_{j,\epsilon}^T=C_{j,\epsilon},$  for $j=1,2,3,$ 

\item  $M_{j,\epsilon}=\left[\ba{rr} D_{j,\epsilon} & B_{j,\epsilon} \\ C_{j,\epsilon} & -D_{j,\epsilon}^T \ea \right]$  is an indecomposable nonsemisimple Hamiltonian \newline $(2N_j \times 2N_j)$-matrix,  for $j=1,2,3$. Moreover, since $N_1,N_2$ are odd and $N_3$ is even,  the matrices $M_{1,\epsilon}, M_{2,\epsilon} $ are given by the formula \eqref{c2} and the matrix $M_{3,\epsilon} $ is given by the formula \eqref{c3}.
\ee
Putting

\be
\item $  \epsilon=-1$ in $M_{1,\epsilon}$,
\item $ \epsilon=+1$ in $M_{2,\epsilon}$,
\item $ \epsilon = + 1,$  in $M_{3,\epsilon}$,
\ee 
we obtain
\be
\item $o_+(\beta_0,J_{2N}H''(0))=0,o_-(\beta_0,J_{2N}H''(0))=1,$
\item $e_+(\beta_0,J_{2N}H''(0))=1,e_-(\beta_0,J_{2N}H''(0))=0.$
\ee
Consequently, we obtain $$o_+(\beta_0,J_{2N}H''(0))-o_-(\beta_0,J_{2N}H''(0))- e_+(\beta_0,J_{2N}H''(0))+e_-(\beta_0,J_{2N}H''(0))=-2.$$

Since $\sign \det H''(0) =1,$ $\ib(0,H') =  1$.  We have just shown that the assumptions of Theorem \ref{gbthoe} are fulfilled.
\eex

%{\blue Nowa wersja tego przykladu 
%\bex
%Let $\beta_0>0$ and $H \in C^2(\bR^{2N},\bR)$ be such that  $J_{2N} H''(0)$ is in the form \eqref{eq:Mj}, where
%\be
%\item $s=3, N=N_1+N_2+N_3$ for $N_1=4k+1,N_2=4k-1, k \in \bN$ and $N_3$ even.
%\item $M_{1,\epsilon}$ is of the form \eqref{c2} with $\epsilon = +1$,
%\item  $M_{2,\epsilon}$ is of the form \eqref{c2} with $\epsilon = -1$,
%\item  $M_{3,\epsilon}$ is of the form \eqref{c3} with $\epsilon = +1$,
%\ee
%For such $H$ we obtain $H'(0)=0, \sigma(J_{2N} H''(0))=\{\pm \beta_0 i\}$  and
%\be
%\item $o_+(\beta_0,J_{2N}H''(0))=0,o_-(\beta_0,J_{2N}H''(0))=1,$
%\item $e_+(\beta_0,J_{2N}H''(0))=1,e_-(\beta_0,J_{2N}H''(0))=0.$
%\ee
%Consequently, 
%$$o_+(\beta_0,J_{2N}H''(0))-o_-(\beta_0,J_{2N}H''(0))- e_+(\beta_0,J_{2N}H''(0))+e_-(\beta_0,J_{2N}H''(0))=-2.$$

%Moreover, since $\det H''(0) =1,$ $\ib(0,H') =  1$.  Therefore the assumptions of Theorem \ref{gbthoe} are fulfilled.
%\eex}

 We emphasize that from the assumptions of  the vast majority of results on   bifurcations of nonstationary periodic solutions of Hamiltonian systems it follows that  there is at least one semisimple eigenvalue of the matrix $J_{2N}H''(0)$ and that the origin is a nondegenerate critical point of the Hamiltonian $H$. We discuss this phenomenon in Propositions \ref{colct} -  \ref{cofrmb}.
 
 Moreover, in these propositions we show that the assumptions of some  well-known theorems imply that of Theorem \ref{thm:main}. It is important to note  here that these assumptions are easy to verify. Here and subsequently $\sign(H''(0))$ stands for the signature of $H''(0)$ i.e.
 $\sign (H''(0))=\morse(-H''(0))-\morse(H''(0))$.

We will consider the following assumptions.

\nt  $\mathbf{(A0)}$ (Assumptions of the
Lyapunov Center Theorem, see Theorem 1.6 of  \cite{Berti}, and also \cite{Lyapunov}) Let $\lambda_1=\beta_0 i, \lambda_2 = -\beta_0 i,\lambda_3,\ldots,\lambda_{2N}$ be the eigenvalues of the matrix $J_{2N}H''(0)$, where $ \beta_0 > 0$.   The following nonresonance condition holds true:
$\ds \frac{\lambda_k}{\lambda_1} \not \in \bZ$ for $k \geq 3$.

\nt  $\mathbf{(A1)}$ (Assumption of the Weinstein theorem, see Theorem 2.1 of  \cite{Weinstein}) The Hessian $ H''(0)$ is positive definite.
 
\nt  $\mathbf{(A2)}$ (Assumptions of the version of the Fadell-Rabinowitz theorem discussed in \cite{Berti}) The Hessian $ H''(0)$ is nondegenerate, the signature $ \sign H''(0)$ is nonzero and there is $T_0 > 0$ such that any nonzero solution of the linearized system $\dot x(t)=J_{2N} H''(0) x(t)$  is   $T_0$-periodic.

To formulate the next two assumptions, we need to represent  $\bR^{2N}$ as a direct sum of two symplectic subspaces.
Assume that    $\bR^{2N}=E_1 \oplus E_2$, where $E_1,E_2$ are symplectic subspaces, invariant  for the flow given by the linearized system $\dot x(t)=J_{2N} H''(0) x(t)$.  Additionally assume that the following condition is fulfilled

\vspace{0.3cm}
\nt $\mathbf{(C)}$ there is $T_0 > 0$ such that all the solutions with initial data in $E_1$ are $T_0$-periodic and there are no $T_0$-periodic solutions in $E_2 \setminus \{0\}$. 

Define a Hamiltonian  $H_i : E_i \to \bR $ as follows 
$\ds  H_i(x_i)= \frac{1}{2}\langle H''(0) x_i,x_i  \rangle$ for $i=1,2$.

\nt  $\mathbf{(A3)}$ (Assumptions of the Moser theorem, see Theorem 4 of \cite{Moser}) The Hessian $H''(0)$ is nondegenerate,  the condition $\mathbf{(C)}$ is fulfilled and the Hessian $H_1''(0)$ is positive definite. 
 
\nt  $\mathbf{(A4)}$ (Assumptions of the Fadell-Rabinowitz theorem, see Theorem 8.4 of \cite{FadRab} and the  Bartsch theorem, see Theorem 1.1 of \cite{Bartsch}) The Hessian $H''(0)$ is nondegenerate,  the condition $\mathbf{(C)}$ is fulfilled and the signature $\sign H_1''(0)$ is nonzero.

\begin{Proposition} \label{colct} If the assumption $\mathbf{(A0)}$ is satisfied,   then $\pm \beta_0 i \in \sigma(J_{2N}H''(0)), \beta_0 >0,$ is a simple eigenvalue and $\ib(0,H')= \pm 1$.  Moreover,   the condition \eqref{eq:sumaLap} is fulfilled.   Hence all the  assumptions of  Theorem \ref{thm:main} are satisfied.
\end{Proposition}
\begin{proof} From the nonresonance condition of the  assumption $\mathbf{(A0)}$ it follows that $\pm \beta_0 i $ is a simple eigenvalue of the matrix $J_{2N}H''(0)$ and $\det H''(0) \not = 0$. It is known that $\det H''(0) \not = 0$ implies $\ib(0,H')=\pm 1$. Since $\pm \beta_0 i $ is  of multiplicity one, Remark \ref{rem:semi} implies that
$e_\pm(\beta_0,J_{2N}H''(0)) =0$ and 

\noindent either $$o_+(\beta_0,J_{2N}H''(0))=1 \text{ and } o_-(\beta_0,J_{2N}H''(0))=0$$ or $$o_+(\beta_0,J_{2N}H''(0))=0 \text{ and } o_-(\beta_0,J_{2N}H''(0))=1.$$ 
   In both cases the condition \eqref{eq:sumaLap} is fulfilled, which completes the proof.

\end{proof}

\begin{Proposition} \label{cow}
If the assumption $\mathbf{(A1)}$ is satisfied, then the spectrum $\sigma(J_{2N} H''(0)) $ consists  only of purely imaginary semisimple eigenvalues and $\ib(0,H')=1$. Moreover,  for arbitrary eigenvalue $\pm \beta_0 i \in \sigma(J_{2N} H''(0)), \beta_0 > 0,$ the condition \eqref{eq:sumaLap} is fulfilled.  Thus all the  assumptions of  Theorem \ref{thm:main} are satisfied.
\end{Proposition}

\begin{proof}
Since the Hessian $H''(0)$ is positive definite, from Theorem 8 in \cite{[HOFZE]} we obtain the existence of the linear symplectic map $L : (\bR^{2N}, \omega) \to (\bR^{2N},\omega)$ and the numbers $0<\beta_1, \ldots, \beta_N$, such that
\begin{equation}\label{eq:positivedef}
\frac{1}{2} \langle (L^T H''(0) L) (x,y)^T, (x,y)^T\rangle  = \frac{1}{2} \sum_{j=1}^N \beta_j(x_j^2+y_j^2).
\end{equation}
Hence the normal form of $J_{2N}H''(0)$ is given by $J_{2N}(L^TH''(0)L).$ Moreover, since $L$ is symplectic,  $\sigma(J_{2N}H''(0))=\sigma(J_{2N}(L^T H''(0)L))$ and therefore $\pm i \beta_j$ are the eigenvalues of  $\sigma(J_{2N}H''(0)).$ Taking into consideration formula \eqref{eq:positivedef} we obtain that all the eigenvalues of $JH''(0)$ are purely imaginary and semisimple and for any eigenvalue  $\pm\beta_{0}i \in \sigma(J_{2N} H''(0))$ with $\beta_0>0$ we have
\begin{enumerate}
\item[1)] $ o_+(\beta_{0},J_{2N}H''(0))=e_+(\beta_{0},J_{2N}H''(0))=e_-(\beta_{0},J_{2N}H''(0))=0,$
\item[2)] $2 \cdot  o_-(\beta_{0},J_{2N}H''(0)) = \mu(\beta_{0})>0,$ where $\mu(\beta_{0})$ is the multiplicity of $\beta_{0}$ considered as an eigenvalue of the matrix $L^T H''(0) L$
\end{enumerate}
This implies condition \eqref{eq:sumaLap}. Additionally, again from formula \eqref{eq:positivedef} we obtain $\ib(0,H')=\sign \det H''(0) =1$. Summing up, all the assumptions of Theorem \ref{thm:main} are fulfilled.

\end{proof}

\begin{Proposition} \label{cofrb}
Assume the Hessian $ H''(0)$ is nondegenerate and  of nonzero   signature $ \sign (H''(0))$. If all the eigenvalues of  the matrix $J_{2N} H''(0)$ are purely imaginary and semisimple, then  there is an eigenvalue  $\pm \beta_0 i\in \sigma(J_{2N} H''(0)),$ $\beta_0 > 0,$ such that  the condition \eqref{eq:sumaLap} is fulfilled. Moreover,  $\ib(0,H')=\pm 1$. Therefore all the  assumptions of  Theorem \ref{thm:main} hold true.
\end{Proposition}
\begin{proof}
First of all notice that  $\ib(0,H') = \sign \det H''(0)  = \pm 1$, since $ H''(0)$ is nondegenerate.

 Denote  the eigenvalues of  the matrix $J_{2N} H''(0)$  by $\pm \beta_1 i,\ldots,\pm \beta_N i,$ where $ \beta_1,\ldots, \beta_N > 0 $, and note that it can happen that $\beta_{j}=\beta_{j'}$ for $j \neq j'$.
 Since all these eigenvalues are purely imaginary and semisimple, there is a symplectic linear map $L : \bR^{2N} \to \bR^{2N}$ such that the Hamiltonian corresponding to the  normal form of $J_{2N}H''(0)$ is given by the following formula

$$\ba{rcl} \frac{1}{2} \langle (L^T H''(0) L) (x,y)^T, (x,y)^T\rangle & = & \ds \frac{1}{2} \sum_{j=1}^p \beta_j(x_j^2+y_j^2) - \frac{1}{2}\sum_{j=p+1}^{N} \beta_j(x_j^2+y_j^2) \ea,$$  see \cite{[AMZE]}, \cite{[HOFZE]}.

Set $\sigma(J_{2N} H''(0))=\{\pm \beta_{j_1}i, \ldots, \pm \beta_{j_s}i\}$.  We claim that $o_+(\beta_{j_0},J_{2N} H''(0)) \not = o_-( \beta_{j_0}, J_{2N} H''(0))$,  for some $j_0 \in \{j_1,\ldots,j_s\}$.

Suppose, contrary to our claim, that  for any   $j_r \in \{j_1, \ldots,j_s\}$ the following equality holds true  
 $o_+(\beta_{j_r},J_{2N} H''(0))  = o_-(\beta_{j_r}, J_{2N} H''(0))$. 
 
 Thus it is easy to see that   for any  $j_r \in\{j_1,\ldots,j_s\}$,
 $$\mu(\beta_{j_r})=2\cdot  o_-(\beta_{j_r} ,J_{2N} H''(0)) = 2 \cdot  o_+(\beta_{j_r} ,J_{2N} H''(0)) =\mu(-\beta_{j_r}),$$ where $\mu(\pm \beta_{j_r})$ is the multiplicity of $\pm \beta_{j_r}$ considered as an eigenvalue of the matrix $L^T H''(0) L$.

Consequently, using the fact that by  Sylvester's law of inertia we have $m^-(H''(0))=m^-(L^T H''(0)L)$, we obtain  $$ \morse(-H''(0))=\sum_{r=1}^s \mu(\beta_{j_r}) = \sum_{r=1}^s \mu(-\beta_{j_r}) =  \morse(H''(0)),$$ which implies $\sign (H''(0))=0,$ a contradiction.

Summing up, there exists $j_0 \in \{j_1,\ldots,j_s\}$ such that 
$$o_+( \beta_{j_0} ,J_{2N} H''(0)) \not = o_-(\beta_{j_0} , J_{2N} H''(0)),$$ $$  e_+(\beta_{j_0} ,J_{2N} H''(0))   = e_-(  \beta_{j_0} , J_{2N} H''(0))=0.$$  
 
Note that we have just shown, that all the assumptions of Theorem \ref{thm:main} are fulfilled for $\beta_0= \beta_{j_0}.$
 \end{proof}

Let us remind that    $\bR^{2N}=E_1 \oplus E_2$, where $E_1,E_2$ are symplectic subspaces, invariant  for the flow given by the linearized system $\dot x(t)=J_{2N} H''(0) x(t)$.
Set $\dim E_i=N_i, i=1,2$ and consider  the linearized Hamiltonian systems $\dot x_i(t)=J_{2N_i}H_i''(0)x_i(t)$ on $\bR^{2N_i}$ for $i=1,2$, where the Hamiltonian  $H_i : E_i \to \bR $ is given by the formula 
$\ds  H_i(x_i)= \frac{1}{2}\langle H''(0) x_i,x_i  \rangle$.

\begin{Proposition} \label{cofrmb}
The following statements  hold true.

\be
\item \label{cofrmb1}
Assume that the Hessian $H''(0)$ is nondegenerate and that the Hessian $H_1''(0)$ is positive definite. If $\sigma (J_{2N_1}H_1''(0)) \cap \sigma(J_{2N_2}H_2''(0)) = \emptyset$, then  for arbitrary  eigenvalue  $\pm \beta_{0} i\in \sigma(J_{2N_1} H_1''(0)),$ $ \beta_0 > 0$,    the condition \eqref{eq:sumaLap} is satisfied.  Moreover $\ib(0,H')=\pm 1$. Therefore all the     assumptions of  Theorem \ref{thm:main} are satisfied.
\item \label{cofrmb2} Assume that the Hessian $H''(0)$ is nondegenerate, and all the eigenvalues of the matrix $J_{2N_1} H_1''(0)$ are purely imaginary and semisimple,  and that  the signature  $\sign H_1''(0)$ is nonzero. If $\sigma (JH_1''(0)) \cap \sigma(JH_2''(0)) = \emptyset$, then    there is an   eigenvalue  $\pm \beta_{0} i\in \sigma(J_{2N_1} H_1''(0)), \beta_0 > 0$,  such that   the condition \eqref{eq:sumaLap} is fulfilled. Additionally, $\ib(0,H')=\pm 1$. That is why all the  assumptions of  Theorem \ref{thm:main} are fulfilled.
\ee
\end{Proposition}
\begin{proof} 
(1)  Since the matrix $H_1''(0)$ is positive definite, repeating the reasoning from the proof of Proposition  \ref{cow} we obtain that    for any eigenvalue  $\beta_{0}i \in \sigma(J_{2N_1} H_1''(0))$ the following conditions hold true
$ o_-(\beta_{0},J_{2N_1}H_1''(0))  >0$ and   
$$ o_+(\beta_{0},J_{2N_1}H'_1(0))=e_+(\beta_{0},J_{2N_1}H_1''(0)) = e_-(\beta_{0},J_{2N_1}H_1''(0))=0.$$   
Finally taking into account that  $\sigma (J_{2N_1}H_1''(0)) \cap \sigma(J_{2N_2}H_2''(0)) = \emptyset$ and  $\sigma (J_{2N_1}H_1''(0)) \cup \sigma(J_{2N_2}H_2''(0)) = \sigma(J_{2N}H''(0))$, we obtain
$$o_\pm(\beta_{0},J_{2N } H''(0))=  o_\pm(\beta_{0},J_{2N_1}H_1''(0)),  e_\pm(\beta_{0},J_{2N } H''(0)) = e_\pm(\beta_{0},J_{2N_1}H_1''(0)).$$    Hence we have proved that condition \eqref{eq:sumaLap} is satisfied. 
Summing up, since $\ib(0,H')  = \sign \det H''(0)= \pm 1$,   all the assumptions of Theorem \ref{thm:main} are fulfilled. \\
(2)   Since  $\sign H''_1(0) \not = 0$,  repeating the reasoning from the proof of Proposition  \ref{cofrb} we obtain that    there is a purely imaginary semisimple eigenvalue  
$\pm \beta_{0}i \in \sigma(J_{2N_1} H_1''(0)),$ $ \beta_0 > 0$, such the following conditions hold true
 
 $o_+(\beta_0,J_{2N_1} H_1''(0)) \not = o_-( \beta_0 , J_{2N_1} H_1''(0)), e_+(\beta_0 ,J_{2N_1} H_1''(0))   = e_-(\beta_0 , J_{2N_1} H_1''(0))=0$.  
 
\nt Finally taking into account that  $$\sigma (J_{2N_1}H_1''(0)) \cap \sigma(J_{2N_2}H_2''(0)) = \emptyset \textrm{ and } \sigma (J_{2N_1}H_1''(0)) \cup \sigma(J_{2N_2}H_2''(0)) = \sigma(J_{2N}H''(0))$$ we obtain
$$ o_\pm(\beta_{0},J_{2N } H''(0))= o_\pm(\beta_{0},J_{2N_1}H_1''(0)),  e_\pm(\beta_{0},J_{2N } H''(0))= e_\pm(\beta_{0},J_{2N_1}H_1''(0)), $$ which completes the proof. Summing up, since $\ib(0,H') = \sign \det H''(0)= \pm 1$,   all the assumptions of Theorem \ref{thm:main} are satisfied.
\end{proof}

 \br \label{general}
Note that
\be
%\item  Corollary \ref{colct} says that the assertion of Theorem \ref{thm:main}  follows from the assumptions of the Lyapunov Center Theorem  $\mathbf{(A0)}$.

%\item  Corollary \ref{cow} says that the assertion of Theorem \ref{thm:main}  follows from the assumptions of the Weinstein theorem  $\mathbf{(A1)}$.

\item The condition that spectrum $\sigma(J_{2N} H''(0))$ consists only of purely imaginary se\-mi\-sim\-ple eigenvalues is equivalent to the condition that all the solutions of the linearized system $\dot x(t)=J_{2N}H''(0) x(t)$ are periodic.Therefore  the assumptions of the Fadell-Rabinowitz theorem  $\mathbf{(A2)}$ are stronger than the assumptions of Proposition \ref{cofrb}, because we do not assume that all the solutions of the system $\dot x(t)=J_{2N}H''(0) x(t)$ are of the same period  $T_0 > 0$. Summing up, from Proposition \ref{cofrb} it follows that the assumption $\mathbf{(A2)}$ implies the assertion of Theorem \ref{thm:main}.

\item In the same way, the assumptions of the Moser theorem $\mathbf{(A3)}$ are stronger than the assumption of Proposition \ref{cofrmb}(\ref{cofrmb1}) and the assumptions of the Fadell-Rabinowitz and Bartsch theorems $\mathbf{(A4)}$ are stronger than the assumption of Proposition \ref{cofrmb}(\ref{cofrmb2}).
\ee
\er

\br
Summing up, under one of the assumptions $\mathbf{(A0)}-\mathbf{(A4)}$, the origin is a nondegenerate critical point of the Hamiltonian $H$ and there exists at least one purely imaginary semisimple eigenvalue $\pm \beta_0 i$ of the matrix  $J_{2N}H''(0)$ such that condition \eqref{eq:sumaLap} is satisfied. Under one of these assumptions, in the papers \cite{Bartsch}, \cite{FadRab}, \cite{Moser}, \cite{Weinstein} the authors have proved a lower estimations of the number of sequences of nonstationary periodic solutions of the system \eqref{has} emanating from the origin.
The important point to note here is that assumptions of these theorems imply that of Theorem \ref{thm:main}. On the other hand, our results are of a different nature. Namely, we study the existence of  closed connected sets  of nonstationary periodic orbits of the system \eqref{has} emanating from the origin.  
\er

\br \label{szulk}
Another result concerning the existence of nonstationary periodic solutions of system \eqref{has} in any neighborhood  of the origin has been proved by Szulkin, see \cite{Szulkin}. The result in the case $\det H''(0) \not = 0$ is given in Theorem 4.1 of this paper. Note that assumptions of this theorem imply the condition \eqref{eq:sumaLap}. Therefore our result generalizes this theorem, since we prove the existence of connected set of solutions, not only the sequence. On the other hand, Theorem 4.4 of \cite{Szulkin} concerns the situationwhen $0$ is an isolated, degenerate critical point of $H$. In this case it is additionally assumed that $c^q(H,0) \not \approx 0$ for some $q$, where $c^q(H,0)$ is a critical group. It can happen, that in such situation $\ib(0,H')$ is trivial. However, if we additionally assume that $\ib(0,H') \not = 0$, then assumptions of Theorem 4.4 imply the condition \eqref{eq:sumaLap}. 
\er

\br 
Note that if $\pm \beta_o i \in \sigma(J_{2N}H''(0))$ is of odd multiplicity then assumptions of Szulkin's theorems are fulfilled, see discussion below Proposition 3.6 of \cite{Szulkin}. Observe that in 
 Example \ref{nonsemisimpleeven} we consider the case of the eigenvalue of even multiplicity. In other words,  Theorem \ref{thm:main} can be considered a generalization of the Lyapunov Center Theorem for Hamiltonian systems, which in some cases is valid for strictly nonsemisimple eigenvalues of even multiplicity and with the origin being isolated degenerate critical point of the Hamiltonian $H$ with nontrivial Brouwer index.
\er

\br
The results concerning the emanation of nonstationary periodic solutions of Hamiltonian systems in the special case of planar and spatial systems with Coriolis forces, have been obtained by direct calculation in \cite{GPRU}. Note that these results can be obtained as a direct consequence of Theorem \ref{thm:main}.
\er

\subsection{Global bifurcation theorem} \label{secgbth}

In this section we formulate sufficient conditions for the existence of  global bifurcations of nonstationary $2\pi$-periodic solutions of the following system
\beq
\label{hamsys}
\dot{x}(t)=\lambda J_{2N}H'(x(t)),
\eeq
where $\Omega \subset \bR^{2N}$ is   open, $H \in C^2(\Omega,\bR)$  and $(H')^{-1}(0) $ is finite. 
In other words, we study properties of continua (i.e. closed connected sets) of $2 \pi$-periodic solutions  of the system \eqref{hamsys}.

The main result of this section is Theorem \ref{gbthoe}. We obtain this theorem as  a direct consequence of Theorem \ref{gbth}. Let us remind that the study of the nontriviality of the bifurcation index   given by  the formula \eqref{bifin} consists of two steps, see the formula \eqref{bifi}. On the one hand, we have to  show the nontriviality of the Brouwer index of the equilibrium point $s_{k_0} $, on the other hand, we have to prove a change of the Morse index of the family of symmetric matrices $T^N_j(\lambda H''(s_{k_0}))$, when the parameter $\lambda$ crosses the value $\ds \frac{j}{\beta_0}$. In  theorem below we express a change of the Morse index of the matrix  $T^N_j(\lambda H''(s_{k_0}))$ in terms of the numbers $o_\pm,e_\pm$ given in Definition \ref{def:numbers}. 
 
\begin{Theorem} \label{gbthoe}
Suppose that the
  Hamiltonian system \eqref{hamsys}  satisfies the conditions $H \in C^2(\Omega, \bR)$ and 
 $(H')^{-1}(0)  =\{s_1, \ldots, s_k\}.$
 
Fix $s_{k_0} \in \{s_1, \ldots, s_k\}$ and assume that 
\be 
\item $\ib(s_{k_0},H') \not = 0$,
\item there exists a pair $\pm i \beta_0$ of purely imaginary eigenvalues of $J_{2N}H''(s_{k_0})$   such that 
\begin{equation}\label{eq:sumaoe}
o_+(\beta_0,J_{2N}H''(s_{k_0}))-o_-(\beta_0,J_{2N}H''(s_{k_0}))-e_+(\beta_0,J_{2N}H''(s_{k_0}))+e_-(\beta_0,J_{2N}H''(s_{k_0})) \neq 0.
\end{equation}
\ee 
Then $(s_{k_0}, \frac{1}{\beta_0}) \in  \cT=\{s_1, \ldots, s_k\} \times (0, +\infty)$ is a global bifurcation point of nonstationary $2\pi$-periodic solutions of the system \eqref{hamsys}.  Moreover, if the continuum $\cC(s_{k_0}, \lambda_0)\subset C_{2\pi}( \Omega) \times (0, +\infty)$ is compact, then
\be 
\item[(a)] $\cC(s_{k_0}, \lambda_0) \cap \cT$ is finite,
\item[(b)] $\displaystyle{\sum_{(\hat{s}, \hat{\lambda})\in \cC(s_{k_0}, \lambda_0) \cap \cT}} \bif(\hat{s},\hat{\lambda})=\Theta.$
\ee
\end{Theorem}

\br 
We remind that the formula \eqref{eq:sumaoe} can be verified directly from the normal form of the matrix $J_{2N}H''(s_{k_0})$. Moreover, it  results from more readable assumptions  $\mathbf{(A0)}$-$\mathbf{(A4)}$ discussed in the previous section.\er

%-----------------------------------------------------------------

\br
The problem of   studying  $2 \pi$-periodic solutions of the family \eqref{hamsys} has  variational nature, see Remark \ref{drv}. Therefore one can prove the existence of a local bifurcation of $2 \pi$-periodic solutions of the system \eqref{hamsys}, see Remark \ref{kb}. To prove Theorem \ref{gbthoe} we have applied the degree for $\sone$-equivariant gradient maps, see Remark \ref{drv}.
The advantage of using the degree for $\sone$-equivariant gradient maps lies in the fact that a change of this degree along the set of stationary solutions  implies the existence of a global bifurcation of $2 \pi$-periodic solutions of the system \eqref{hamsys}.
The choice of the degree for $\sone$-equivariant gradient maps seems to be the best adapted to our theory.
\er

%------------------------------------------------------------------

\section{Proof of the main result}\label{sec:proof}

In this section we prove the main results of our article.

\subsection{Proof of Theorem \ref{gbthoe}}
Let $s_{k_0} \in (H')^{-1}(0) $ satisfy assumptions (1) and (2) of Theorem \ref{gbthoe} and fix $\pm \beta_0 i \in \sigma(J_{2N}H''(s_{k_0})), \beta_0>0$ such that condition \eqref{eq:sumaoe} holds. We claim that  $\bif(s_{k_0}, \frac{1}{\beta_0}) \neq \Theta,$ where $\bif(s_{k_0}, \frac{1}{\beta_0})$ is defined by formula \eqref{bifin} and therefore, our assertion follows from Theorem \ref{gbth}.

Let $T_1^N(\lambda H''(s_{k_0}))$ be given by \eqref{eq:Tj}. Define 

\beq \label{jump} 
\gamma \left( \frac{1}{\beta_0}\right) =  \morse \left(T_1^N \left( \left(\frac{1}{\beta_0} + \mu \right)H''(s_{k_0}) \right) \right) - \morse \left(T_1^N \left( \left(\frac{1}{\beta_0} - \mu \right) H''(s_{k_0}) \right) \right), 
\eeq
where
$\mu$ is chosen such that $\ds \left[\frac{1}{\beta_0}-\mu,\frac{1}{\beta_0}+\mu \right] \cap \Lambda(s_{k_0}) = \left\{\frac{1}{\beta_0}\right\}.$ 

Since, by assumption, $\ib(s_{k_0}, H') \neq 0,$ to complete  the proof it is enough to show that $\gamma \left( \frac{1}{\beta_0}\right)\neq 0$.
Therefore in the rest of  this section we  compute a change of the Morse index of the matrix $T_1^N(\lambda H''(s_{k_0}))$ with respect to a positive parameter $\lambda.$

\begin{Remark}\label{rem:Tk} Note that if Hamiltonian matrices $M_1=J_{2N} A_1$ and $M_2=J_{2N} A_2$ are symplectically similar, then  $\morse(T_1^N(\lambda A_1))=\morse(T_1^N(\lambda A_2))$.
Indeed, if $J_{2N}A_2=S^{-1}J_{2N}A_1S$, then $A_2=S^TA_1 S$ and therefore $T_1^N(\lambda A_2)=\bS^T T_1^N(\lambda A_1) \bS$, where $\bS=\left[\begin{array}{cc} S &0\\ 0 &S\end{array}\right].$ Therefore the assertion follows from the Sylvester's law of inertia.
\end{Remark}

To shorten notation, put $A=H''(s_{k_0})$. Taking into account the above remark, from now on without loss of generality we assume that $J_{2N}A$ is in the normal form. 

The following fact is well-known, see for example \cite{BSz} (or Corollary 11.16 of \cite{GPRU}).

\begin{Fact}\label{degeneracy}
The matrix $T_1^N(\lambda A)$ is degenerate iff $\ds \lambda = \frac{1}{\beta}$, where $\pm i \beta  \in \sigma(J_{2N} A).$
\end{Fact}

Our aim is to study a change of the Morse index $\morse (T_1^N(\lambda A))$  when the parameter $\lambda$ crosses the value $\ds \frac{1}{\beta_0}.$

\bl \label{morsesum} Let $A$ be such that $M=J_{2N}A$ is in the form \eqref{eq:Mj} of Remark \ref{rem:normalform} and let $N_j$ and $M_j$, for $j=1, \ldots, s$ be given by this remark. Put $A_j=J_{2N_j}^{-1}M_j$ and $T_1^{N_j}(\lambda A_j)=\left[\ba{cccc} -\lambda A_j& J_{2N_j} \\ -J_{2N_j} & -\lambda A_j \ea \right].$  Then
 \beq \label{morsedec} \morse\left(T_1^N\left(\left( \frac{1}{\beta_0 }\pm \mu\right) A \right)\right) = \sum_{j=1}^s  \morse\left(T_1^{N_j} \left(\left( \frac{1}{\beta_0} \pm \mu\right) A_j\right)\right).\eeq 
\el

\begin{proof}
Set $x=(x_1, \ldots, x_{2s}),y=(y_1,\ldots,x_{2s}) \in \bR^{2N}=\bR^{N_1} \oplus \ldots \oplus \bR^{N_s} \oplus \bR^{N_1} \oplus \ldots \oplus \bR^{N_s}$ i.e. $x_i, x_{i+s},y_i,y_{i+s} \in \bR^{N_i}$ for $i=1,\ldots,s.$
Let us define $\bV_j=\{(x,y) \in \bR^{4N} : x_i = 0, y_i=0 \text{ for } i \not = j, j+ s\}$ , for $j=1, \ldots,s.$

We obtain
\be
\item $\bR^{4N}=\bV_1 \oplus \ldots \oplus \bV_s,$
\item for $j=1, \ldots,s$,  $\bV_j$ is an invariant subspace of $T_1^N\left(\lambda A \right).$
\ee
Therefore, there exist an orthogonal change of coordinates transforming $T_1^N(\lambda A)$ into the matrix $\diag\left( T_1^{N_1} (\lambda A_1), \ldots, T_1^{N_s}(\lambda A_s)\right)$, which ends the proof.

\end{proof}

\br \label{decomp}  For $j=1, \ldots,s$ define $$ \gamma_j \left(\frac{1}{\beta_0}\right)=   \morse \left(T_1^{N_j} \left( \left(\frac{1}{\beta_0} + \mu \right) A_j \right) \right) - \morse \left(T_1^{N_j} \left( \left(\frac{1}{\beta_0} - \mu \right)A_j \right) \right),$$ where $\mu$ is chosen as in formula \eqref{jump}. By  Lemma \ref{morsesum} and the formula  \eqref{jump} we obtain
\beq \label{sumjump} \gamma\left(\frac{1}{\beta_0} \right) = \sum_{j=1}^s \gamma_j \left(\frac{1}{\beta_0}\right).
\eeq
\er

In view of the above remark, from now on we fix $j \in \{1, \ldots, s\}$ and consider $A_j$ given by Lemma \ref{morsesum}.

\br \label{when0}
From Fact \ref{degeneracy} it follows that if $\pm \beta_0 i  \notin \sigma(J_{2N_j} A_j)$, then
$\ds \gamma_j\left(\frac{1}{\beta_0}\right)=0.$  Therefore, since  $M_j=J_{2N_j}A_j$ is an indecomposable Hamiltonian matrix, what is left is to compute $\ds \gamma_j\left(\frac{1}{\beta_0}\right)$ for matrices $A_j$ satisfying $\sigma(J_{2N_j} A_j)=\{\pm \beta_0 i\}$. 
\er

\bl \label{morseindex} If $\sigma(J_{2N_j} A_j)=\{\pm \beta_0 i\}, \beta_0 > 0,$ then 
$$\morse \left( T_1^{N_j}\left( \left( \frac{1}{\beta_0}+\mu \right) A_j \right) \right) = \left\{\ba{ccc} 2N_j & \text{ for } & \mu < 0  \\ 2 \morse(-A_j) & \text{ for } & \mu > 0\ea  \right..$$    \el
\begin{proof}
Combining Fact \ref{degeneracy} and the assumptions we obtain that  $\det   T_1^{N_j}\left( \left( \frac{1}{\beta_0}+\mu \right) A_j \right)=0$ iff  $\mu=0.$ Hence a change of the Morse index of $ T_1^{N_j}\left( \left( \frac{1}{\beta_0}+\mu \right) A_j \right)$ can occur only at  $\mu = 0.$  That is why for $\mu < 0$ the Morse index  $\morse \left( T_1^{N_j}\left( \left( \frac{1}{\beta_0}+\mu \right) A_j \right) \right)$ does not depend on the choice of $\mu$ and its value can be obtained for example for $\mu=-\frac{1}{\beta_0}$. In this case we obtain $\morse(T^{N_j}_1(\Theta))=2N_j$, where $\Theta$ is the zero $(2N_j \times 2N_j)$-matrix. 

Suppose now that $\mu > 0$ and note that $$\morse \left(  T_1^{N_j}\left( \left( \frac{1}{\beta_0}+\mu \right) A_j \right) \right) = \morse \left(\left( \frac{1}{\beta_0}+\mu \right)^{-1} \cdot T_1^{N_j}\left( \left( \frac{1}{\beta_0}+\mu \right) A_j \right)  \right) =$$ $$= \morse \left(\left[ \ba{cc} -A_j &  \left( \frac{1}{\beta_0}+\mu \right)^{-1}J_{2N_j}\\  -\left( \frac{1}{\beta_0}+\mu \right)^{-1}J_{2N_j}& -A_j  \ea \right]\right) .$$ Since the Morse index does not depend on the choice of positive $\mu$ and the matrix obtained as the limit for $\mu \to \infty$ is nondegenerate, we obtain 
$$\morse \left(\left[ \ba{cc} -A_j &  \left( \frac{1}{\beta_0}+\mu \right)^{-1}J_{2N_j}\\  -\left( \frac{1}{\beta_0}+\mu \right)^{-1}J_{2N_j}& -A_j  \ea \right]\right)=\morse \left(\left[ \ba{cc} -A_j & \Theta\\  \Theta& -A_j  \ea \right]\right)=2\morse(-A_j)$$
 which completes the proof.
\end{proof}

As a direct consequence of the above lemma we obtain the following corollary.

 \bco  \label{jumpA}
If $\sigma(J_{2N_j} A_j)=\{\pm \beta_0 i\}, \beta_0 > 0$, then $$\gamma_j \left( \frac{1}{\beta_0} \right)=2 (\morse(-A_j) -N_j).$$
\eco

From the above considerations it follows that to compute $\gamma \left( \frac{1}{\beta_0} \right)$ given by \eqref{jump}, it is enough to compute the numbers $\gamma_j \left( \frac{1}{\beta_0} \right)$ defined in Remark \ref{decomp}.  Moreover, taking into account Remark \ref{when0} it follows that one can restrict the reasoning to the cases when $ \{\pm \beta_0 i\} = \sigma(J_{2N_j}A_j)$.  Since the considered Hamiltonian matrices are indecomposable, we only have to study cases (C1),(C2) given in Section \ref{subs:normal}.  Without loss of generality, in two following lemmas we assume that $s_{k_0}=0.$

%\bl \label{simple} In the case (C1),  $\gamma \left(\frac{1}{\beta_0} \right) = -2 \epsilon.$  \el
%\begin{proof}
%Since $\ds -H''(0)= \left[\ba{rr}  \epsilon \beta_0 & 0 \\ 0 & \epsilon \beta_0  \ea \right]$ By Corollary \ref{jumpA} we obtain   $$\gamma \left( \frac{1}{\beta_0 } \right) = 2 (\morse(-H''(0)) -N) = 2 ((1-\epsilon) - 1)=-2 \epsilon,$$ which completes the proof.
%\end{proof}

\bl \label{nonsimpleO} Let $\epsilon \in \{+1,-1\}$ and  $A_j=J_{2N_j}^{-1}M_{j,\epsilon}$, where $M_{j,\epsilon}$ is given  in the case (C1). Then  $\gamma_j \left(\frac{1}{\beta_0} \right) = 2 (-1)^{\frac{N_j+1}{2}}\epsilon.$  \el
\begin{proof}
Define a family of Hamiltonians $\cH_{\tau}(x_1,\ldots,x_{N_j},y_1,\ldots,y_{N_j}),$ $ \tau \in [0,1],$ by the following formula

 $$\cH_{\tau}(x_1,\ldots,x_{N_j},y_1,\ldots,y_{N_j})=$$ $$= -\epsilon \beta_0  \sum_{i=1}^{\frac{N_j-1}{2}} (-1)^{i+1} (x_i x_{N_j+1-i} + y_i y_{N_j+1-i}) + \tau \sum_{i=1}^{N_j-1} x_i y_{i+1} + \frac{\epsilon \beta_0}{2} (-1)^{\left[\frac{N_j}{2} \right] +1} \left(x^2_{\frac{N_j+1}{2}} + y^2_{\frac{N_j+1}{2}}\right).$$

By elementary computations and reasoning by induction it is easy to verify that  
\be
\item \label{h1h} $\cH_1(x,y)=\frac{1}{2}\langle A_j(x,y)^T, (x,y)^T \rangle,$ where $(x,y)=(x_1,\ldots,x_{N_j},y_1,\ldots,y_{N_j})$,
\item \label{dettau} $\det \cH''_{\tau}(0)= \beta_0^{2N_j}$ for all $\tau \in [0,1]$,
\item \label{deter} $\det (-\cH_0''(0)-\lambda Id)=\left( \lambda+\epsilon \beta_0 \right) ^{N_j-1} \left( \lambda - \epsilon \beta_0 \right) ^{N_j-1}(\lambda-(-1)^{\frac{N_j-1}{2}}\epsilon \beta_0)^2.$
\ee

By \eqref{dettau}  $-\cH''_{\tau}(0)$ is a path of isomorphisms and that is why the Morse index $\morse (-\cH_{\tau}''(0))$ does not change along this path.
Hence by \eqref{h1h} we obtain $$\morse (-A_j)= \morse (-\cH_1''(0))= \morse (-\cH_0''(0)).$$
In other words we have simplified the computations of the Morse index $\morse (-A_j).$

 Taking into account \eqref{deter} we obtain that  $$ \morse(-A_j)=N_j-1+(1-(-1)^{\frac{N_j-1}{2}} \epsilon)= N_j + (-1)^{\frac{N+1}{2}} \epsilon.$$

Finally by Corollary \ref{jumpA} we obtain
$$\gamma_j \left(\frac{1}{\beta_0}\right)=2(\morse(-A_j)-N_j)=2 (-1)^{\frac{N_j+1}{2}} \epsilon,$$
 which completes the proof.
\end{proof} 

\bl \label{nonsimpleE}  Let $\epsilon \in \{+1,-1\}$ and $A_j=J_{2N_j}^{-1}M_{j,\epsilon}$, where $M_{j,\epsilon}$ is given  in the case (C2). Then  $\gamma_j \left(\frac{1}{\beta_0} \right) = 0.$  \el
\begin{proof}
In view of Corollary \ref{jumpA}  it is enough to show that  $\morse(-A_j)=N_j.$ Define a family of Hamiltonians $\cH_{\tau}(x_1,\ldots,x_{N_j},y_1,\ldots,y_{N_j}), \tau \in [0,1],$ by the following formula
 $$\cH_{\tau}(x_1,\ldots,x_{N_j},y_1,\ldots,y_{N_j})=\beta_0 \sum_{i=1}^{N_j/2} (x_{2i-1}y_{2i} - x_{2i}y_{2i-1}) + \tau \left(\sum_{i=1}^{N_j-2} x_i y_{i+2} + \frac{1}{2} \epsilon (x_{N_j-1}^2 + x_{N_j}^2)\right).$$

Reasoning by induction it is easy to verify that
\be
\item $\cH_1(x,y)=\frac{1}{2}\langle A_j(x,y)^T, (x,y)^T \rangle,$ where $(x,y)=(x_1,\ldots,x_{N_j},y_1,\ldots,y_{N_j})$,
\item  $\det \cH''_{\tau}(0)= \beta_0^{2N_j}$ for all $\tau \in [0,1]$,
\item \label{deter1} $\det (-\cH_0''(0)-\lambda Id)=\left( \lambda+ \beta_0 \right) ^{N_j} \left( \lambda - \beta_0 \right)^{N_j}.$
\ee
Repeating the reasoning used in the proof of the previous lemma we obtain
$$\ds \gamma_j \left(\frac{1}{\beta_0}\right)=2(\morse(-A_j)-N_j)=0,$$
 which completes the proof.
\end{proof}

Now we are in the position to complete the proof of Theorem \ref{gbthoe}. To this end we show that assumptions of Theorem \ref{gbthoe} imply that $\bif(s_{k_0}, \frac{1}{\beta_0})$ is nontrivial.

Consider $\eta_1(s_{k_0}, \frac{1}{\beta_0})$, i.e. the first coordinate of $\bif(s_{k_0}, \frac{1}{\beta_0})$. Observe that $\eta_1(s_{k_0}, \frac{1}{\beta_0}) = \ib(s_{k_0}, H') \cdot \gamma(\frac{1}{\beta_0})$. Since $\ib(s_{k_0}, H')\neq 0$, to show nontriviality of $\bif(s_{k_0}, \frac{1}{\beta_0})$ it is enough to prove that $\gamma(\frac{1}{\beta_0})\neq 0.$

From Remark \ref{decomp} it follows that $\ds \gamma\left(\frac{1}{\beta_0}\right)=\sum_{j=1}^{s} \gamma_j\left(\frac{1}{\beta_0}\right)$. Moreover, Lemmas \ref{nonsimpleO} and \ref{nonsimpleE} provide formulas for numbers $\gamma_j\left(\frac{1}{\beta_0}\right)$ for each $j \in \{1, \ldots, s\}.$ In particular, these lemmas allow us to split all the matrices $A_j$ in the normal form \eqref{eq:Mj} of $A$ into five groups:
\begin{enumerate}[(1)]
\item matrices with $N_j$ odd, $\frac{N_j+1}{2}$ odd and $\epsilon = +1$. In this case, by Lemma \ref{nonsimpleO}, $\gamma_j\left(\frac{1}{\beta_0}\right)=-2,$
\item matrices with $N_j$ odd, $\frac{N_j+1}{2}$ odd and $\epsilon = -1$. In this case, by Lemma \ref{nonsimpleO}, $\gamma_j\left(\frac{1}{\beta_0}\right)=2,$
\item matrices with $N_j$ odd, $\frac{N_j+1}{2}$ even and $\epsilon = +1$. In this case, by Lemma \ref{nonsimpleO}, $\gamma_j\left(\frac{1}{\beta_0}\right)=2,$
\item matrices with $N_j$ odd, $\frac{N_j+1}{2}$ even and $\epsilon = -1$. In this case, by Lemma \ref{nonsimpleO}, $\gamma_j\left(\frac{1}{\beta_0}\right)=-2,$
\item matrices with $N_j$ even. In this case, by Lemma \ref{nonsimpleE}, $\gamma_j\left(\frac{1}{\beta_0}\right)=0.$
\end{enumerate}

By the definition of numbers $o_{\pm}(\beta_0,J_{2N}H''(s_{k_0})), e_{\pm}(\beta_0,J_{2N}H''(s_{k_0}))$, the number of matrices in the groups (1)-(4) equals, respectively, $o_+(\beta_0,J_{2N}H''(s_{k_0})), o_-(\beta_0,J_{2N}H''(s_{k_0}))$, $e_+(\beta_0,J_{2N}H''(s_{k_0})), e_-(\beta_0,J_{2N}H''(s_{k_0})).$

Finally, we obtain 
\begin{equation*}
\begin{split}
&\gamma \left(\frac{1}{\beta_0} \right) = \sum_{j=1}^s \gamma_j \left(\frac{1}{\beta_0}\right)= \\ =-2(o_+(\beta_0,J_{2N}H''(s_{k_0}))-o_-(\beta_0,&J_{2N}H''(s_{k_0}))-e_+(\beta_0,J_{2N}H''(s_{k_0}))+e_-(\beta_0,J_{2N}H''(s_{k_0}))).
\end{split}
\end{equation*}

Hence, condition \eqref{eq:sumaoe} implies $\bif(s_{k_0}, \frac{1}{\beta_0}) \neq \Theta.$ Therefore our assertion follows from Theorem \ref{gbth}.

\subsection{Proof of Theorem \ref{thm:main}}
Without loss of generality we assume that $(H')^{-1}(0) \cap \Omega =\{0\}.$ Moreover, using the standard change of variables $\tilde{u}(t)=u(\lambda t)$ we can change the problem of studying $2\pi \lambda$-periodic solution  of the system \eqref{has}  to the study of $2\pi$-periodic solutions of the family \eqref{hamsys}. In this case it is obvious that we can restrict our considerations to $\lambda \in (0, +\infty).$

Note that if the assumptions of Theorem \ref{thm:main} are satisfied, then also assumptions of Theorem \ref{gbthoe} are fulfilled. Therefore we obtain a closed, connected set of solutions of system \eqref{hamsys}, bifurcating from $(0, \frac{1}{\beta_0}).$ It remains to show that this implies the existence of the connected set of closed trajectories, emanating from the origin, and such that condition \eqref{crezo} is satisfied. Considering the continous mapping $f \colon C_{2\pi}(\Omega) \to \bR^{2N}$ given by $f(u)=u(0)$ and taking the sufficiently small neighbourhood $\cU$ of $(0, \lambda_0)$,  we obtain the connected set $f(\cC(0, \frac{1}{\beta_0})\cap \cU) \subset \bR^{2N}.$ This implies that the set $\{u([0, 2\pi]); u \in \cC(0, \frac{1}{\beta_0})\cap \cU\}$ is also connected. It is easy to observe, that for such a set the condition \eqref{crezo} is satisfied, which ends the proof.

\vspace{0.4cm}
\noindent \textbf{Data availability statement.} \\
Data sharing not applicable to this article as no datasets were generated or analysed during the current study.

\end{document}